\documentclass[11pt]{article}
\usepackage{amsfonts,latexsym,rawfonts,amsmath,amssymb,amsthm,graphicx}
\numberwithin{equation}{section}
\newtheorem{theorem}{Theorem}[section]

\newtheorem{thm}[theorem]{Theorem}

\newtheorem{defi}[theorem]{Definition}

\def\s{\,\,\,\,}

\def\endproof{$\hfill\Box$\\}

\title{ Weakly convex biharmonic hypersurfaces in nonpositive curvature space forms are minimal}
\author{Yong Luo\footnote{The author is supported by the DFG Collaborative Research Center SFB/Transregio 71.}}
\date{}
\begin{document}
\maketitle
\begin{abstract}
A submanifold $M^m$ of a Euclidean space $R^{m+p}$ is said to have harmonic mean curvature vector field if $\Delta \vec{H}=0$, where $\vec{H}$ is the mean curvature vector field of $M\hookrightarrow R^{m+p}$ and $\Delta$ is the rough Laplacian on $M$. There is a conjecture named after Bangyen Chen which states that submanifolds of Euclidean spaces with harmonic mean curvature vector fields are minimal. In this paper we prove that weakly convex hypersurfaces (i.e. hypersurfaces whose principle curvatures are nonnegative)  with  harmonic mean curvature vector fields in Euclidean spaces are minimal. Furthermore we prove that weakly convex biharmonic hypersurfaces in nonpositive curved space forms are minimal.
\end{abstract}

\section{Introduction}
Let $\vec{x}: M^m\to R^{m+p}$ be an immersion from a Riemmanian manifold $M$ of dimension $m$ to a Euclidean space of dimension $m+p$, $p\geq1$. Denote by $\vec{x}, \vec{H}, \Delta$ respectively the position vector of $M$, the mean curvature vector field of $M$ and the Laplacian operator  with respect to the induced metric $g$ on $M$. Then it is well known that (see for example \cite{Ch1})
$$\Delta \vec{x}=-n\vec{H}.$$
This shows that $M$ is a minimal submanifold if and only if its coordinates functions are harmonic functions. According to this equation, a submanifold in a Euclidean space with harmonic mean curvature vector field, i.e.
\begin{eqnarray}
\Delta \vec{H}=0,
\end{eqnarray}
if and only if
\begin{eqnarray}
\Delta^2 \vec{x}=0.
\end{eqnarray}
Therefore a submanifold with harmonic mean curvature vector field is called a biharmonic submanifold.

There is also a variational description of biharmonic submanifolds as follows. Assume that $\vec{x}: M\to N^{m+p}$ is an immersion to a Riemnnian manifold, then the biharmonic energy of $\vec{x}$ is defined by
\begin{eqnarray}
E_2(\vec{x})=\int_M|\tau(\vec{x})|^2dM,
\end{eqnarray}
where $\tau(\vec{x})$ is the tension field of $\vec{x}$. The critical points of the functional $E_2$ satisfy the following E-L equation (see \cite{Ji})
\begin{eqnarray}\label{equ3}
-\Delta \vec{H}=\sum_{i=1}^mR^N(e_i,\vec{H})e_i,
\end{eqnarray}
where $\{e_i\}$ is a local orthonormal frame on $M$. In particular, when $N$ is a Euclidean space, this equation is just $\Delta \vec{H}=0$.
\begin{defi}
A Submanifold  satisfying equation (\ref{equ3}) is called a biharmonic submanifold.
\end{defi}
As is easy to see, minimal submanifolds are biharmonic submanifolds.
It is nature to arise the question whether the space of biharmonic submanifolds is strictly larger than the space of minimal submanifolds. Concerning this problem, when the ambient manifold is a Euclidean space Chen conjectured that the answer is no.
\\\textbf{Chen's Conjecture.} Suppose that $\vec{x}: M^m\to R^{m+p}$, $p\geq1$, satisfies
\begin{eqnarray}
\Delta \vec{H}=0.
\end{eqnarray}
Then $\vec{H}=0$.

Since Chen's conjecture arose, it remains to be open with very little progress even for hypersurfaces with dimensions grater that 4. But in recent years it attracts
many attentions and there are some partial answers to this conjecture.
Now we give an overview of them, as to the author's knowledge. Chen's conjecture is proved:

$\bullet$ For Surfaces in $R^3$ in \cite{Ch-Is} and \cite{Ji} and in an unpublished work by Chen himself, as reported in \cite{Ch}.

$\bullet$ For hypersurfaces in $R^4$  in \cite{Ha-Vl} and a different proof in \cite{De}.

$\bullet$ For hypersurfaces which admit at most 2 distinct principle curvatures in \cite{Di}.

$\bullet$ For curves \cite{Di}.

$\bullet$ For submanifolds $M^m$ which are pseudo-umbilic and $m\neq 4$ \cite{Di}.

$\bullet$ For submanifolds of finite type \cite{Di}.

$\bullet$  For submanifolds which are proper, i.e. any preimage of compact subsets are compact in \cite{Ak-Ma}.

In this paper we prove Chen's conjecture for weakly convex biharmonic hypersurfaces in Euclidean spaces.
\begin{thm}\label{thm}
Assume that $\vec{x}: M^m\to R^{m+1}$ is a weakly convex biharmonic
submanifold in $R^{m+1}$, i.e. $\Delta \vec{H}=0.$ Then $\vec{H}=0$.
\end{thm}
Chen's conjecture is generalized to the so called generalized Chen's conjecture (see \cite{Ca-Mo} \cite{Ch} \cite{Ch1} \cite{Ou} \cite{Ou-Ta} and \cite{On} etc.) allowing the ambient manifolds to be nonpositive curved.
\\\textbf{Generalized Chen's Conjecture.} Suppose that $f: M\to (N, h)$ is an isometric immersion and the section curvature of $(N, h)$  nonpositive. Then if $f$ is biharmonic, it is minimal.

The generalized Chen's curvature turned out to be false by a counter example constructed by Y-L. Ou and L. Tang (see \cite{Ou-Ta}). But it remains interesting to find out sufficient conditions which make it be true. It is easy to see from the proof of theorem \ref{thm} that our argument is also applicable to find out that for hypersurfaces weakly convex is a sufficient condition to guarantee the generalized Chen's conjecture to be true when the ambient manifold is a space form $N(c)$ of constant curvature $c\leq0$. Thus we have
\begin{thm}\label{thm2}
Assume that $\vec{x}: M^m\to N^{m+1}(c)$ is a weakly convex biharmonic
hypersurfaces in a spcae form $N^{m+1}(c)$ with $c\leq 0$, then $\vec{H}=0$.
\end{thm}
Other sufficient conditions have also been found out to guarantee the generalized Chen's conjecture to be true. For recent development in this direction, we refer to \cite{Na-Ur} \cite{Luo} etc..

We would like to point out that  for closed (compact without boundary)
bi-harmonic submanifolds the Chen's conjecture and the generalized Chen's conjecture  is easily proved to
be true by an argument due to Jiang (see \cite{Ji}) or by directly using
an integration by parts argument.

The Chen's and generalized Chen's conjectures belong to the research field of classification of biharmonic submanifolds in Riemannian or pseudo-Riemannian manifolds. Nowadays it is an active research field and for readers who have interest with it we refer to a survey paper \cite{MO} by Montaldo and Oniciuc and references therein.

\textbf{Organization.} In section 2 we give a brief sketch of submanifolds geometry.
Theorem \ref{thm} and theorem \ref{thm2} are proved in section 3.
\section{Preliminaries}
Assume that $\vec{x}:M^m\to N^{m+p}$ is an immersion to a Riemannian manifold $N$
with Riemannian metric $\langle,\rangle$, which is a bilinear form on $TN\otimes TN$,
the tensor of the tangent vector space of $N$. Then $M$ inherits a
Riemannian metric from $N$ by $g_{ij}:=\langle\partial_i \vec{x}, \partial_j\vec{x}\rangle$
and a volume form by $\sqrt{\det g_{ij}}dx$.
The second fundamental form of $M\hookrightarrow N$, $ h: TM\otimes TM\to NM$, is defined by
$$h(X, Y):=D_XY-\nabla_XY,$$
for any $X, Y\in TM$, where $D$ is the covariant derivative with respect to the Levi-Civita connection on $N$, $\nabla$ is the Levi-Civita connection on $M$ with respect to the induced metric and $NM$ is the normal bundle of $M$. For any normal vector field $\eta$ the Weingarten map $A_\eta: TM\to TM$ is defined by
\begin{eqnarray}
D_X\eta=-A_\eta X+\nabla^\bot_X\eta,
\end{eqnarray}
where $\nabla^\bot$ is the normal connection and as is well known that $h$ and $A$ are related by
\begin{eqnarray}
\langle h(X, Y), \eta\rangle=\langle A_\eta X, Y\rangle.
\end{eqnarray}

For any $p\in M$, let $\{e_1, e_2,...,e_m, e_{m+1},...,e_{m+p}\}$ be a local orthonormal basis of $N$ such that $\{e_1,...,e_m\}$ is an orthonormal basis of $T_pM$. Then $h$ is decomposed  at $p$ as
$$h(X,Y)=\sum_{\alpha=m+1}^{m+p}h_\alpha(X, Y)e_\alpha.$$
The mean curvature vector field is defined as
\begin{eqnarray}
\vec{H}:=\frac{1}{m}\sum_{i=1}^mh(e_i, e_i)=\sum_{\alpha=m+1}^{m+p}H_\alpha e_\alpha,
\end{eqnarray}
where
\begin{eqnarray}
H_\alpha:=\frac{1}{m}\sum_{i=1}^mh_\alpha(e_i, e_i).
\end{eqnarray}
\section{Proof of theorem \ref{thm}}
  It is well known that for a submanifold $M^m$ in a Euclidean space to be biharmonic, i.e. $\Delta \vec{H}=0$ if and only if (\cite{Ch})
  \begin{eqnarray}
\Delta^\bot \vec{H}-\sum_{i=1}^mh(A_{\vec{H}}e_i, e_i)&=&0, \label{equ1}
\\m\nabla|\vec{H}|^2+4\sum_{i=1}^mA_{\nabla^\bot_{e_i} \vec{H}}e_i&=&0, \label{equ2}
\end{eqnarray}
where $\Delta^\bot$ is the (nonpositive) Laplace operator associated with the normal connection $\nabla^\bot$.

 Assume that $\vec{H}=H\nu$, where $\nu$ is the unit normal vector field on $M$. Note that by the assumption that $M$ is weakly convex, we have $H\geq0$. Define
 \begin{eqnarray}
 B:=\{p\in M: H(p)> 0\}
 \end{eqnarray}
We will prove that $B$ is an empty set by a contradiction argument, and so $M$ is minimal and we are done.

If $B$ is not empty, we see that $B$ is an open subset of $M$.  We assume that $B_1$ is a nonempty connect component of $B$. We will prove that $H\equiv 0$ in $B_1$, thus a contradiction.

We prove it in two steps.

 \textbf{Step 1.} $H$ is a constant in $B_1$.

 Let $p\in B_1$ be a point. Around $p$ we choose a local orthonormal frame $\{e_k, k=1,...,m\}$ such that $\langle h, \nu\rangle$ is a diagonal matrix at $p$, where $\nu$ is the unit normal vector field of $M$.

 For any $1\leq k\leq m$ we have at $p$
\begin{eqnarray*}
&&\langle\sum_{i=1}^mA_{\nabla^\bot_{e_i} \vec{H}}e_i, e_k\rangle
\\&=&\sum_{i=1}^m\langle h(e_i, e_k), \nabla^\bot_{e_i}\vec{H} \rangle
\\&=&\langle h(e_k, e_k), \nabla^\bot_{e_k} \vec{H} \rangle.
\end{eqnarray*}
By equation (\ref{equ2}) we have
\begin{eqnarray*}
0&=&m\nabla_{e_k}|\vec{H}|^2+4\langle\sum_{i=1}^mA_{\nabla^\bot_{e_i}\vec{H}}e_i, e_k\rangle
\\&=&m\nabla_{e_k}|\vec{H}|^2+4\langle h(e_k, e_k), \nabla^\bot_{e_k} \vec{H} \rangle
\\&=&2mH\nabla_{e_k}H+4\lambda_k\langle\nu, \nabla^\bot_{e_k}\vec{H}  \rangle
\\&=&(2mH+4\lambda_k)\nabla_{e_k}H,
\end{eqnarray*}
where $\lambda_k:=h(e_k, e_k)$ is the $k$th principle curvature of $M$ at $p$, which is nonnegative by the assumption that $M$ is weakly convex.

By $2mH+4\lambda_k>0$ at $p$, one gets
\begin{eqnarray}
\nabla_{e_k}H=0\s at\s p,
\end{eqnarray}
for any $k=1,...,m$, which implies that
\begin{eqnarray}
\nabla H=0 \s at \s p.
\end{eqnarray}
Because $p$ is an arbitrary point in $B_1$, we see that
\begin{eqnarray}
\nabla H=0 \s in\s B_1.
\end{eqnarray}
Therefore we get that $H$ is a constant in $B_1$.

\textbf{Step 2.} $H$ is zero in $B_1$.

Let $p\in B_1$, by step 1, we see that
\begin{eqnarray}\label{eq1}
\Delta |\vec{H}|^2(p)=0.
\end{eqnarray}
On the other hand, by equation (\ref{equ1}), we have
\begin{eqnarray}\label{ine2}
\Delta |\vec{H}|^2&=&2|\nabla^\bot \vec{H}|^2+2\langle \vec{H}, \Delta^\bot \vec{H}\rangle \nonumber
\\&\geq&2\sum_{i=1}^m\langle h(A_{\vec{H}}e_i, e_i), \vec{H}\rangle \nonumber
\\&=&2\sum_{i=1}^m \langle A_{\vec{H}}e_i, A_{\vec{H}}e_i\rangle \nonumber
\\&=&2\sum_{i=1}^m H^2 \langle A_\nu e_i, A_\nu e_i\rangle \nonumber
\\&\geq&2n H^4.
\end{eqnarray}
From (\ref{eq1})-(\ref{ine2}), we get $H(p)=0$. Because $p$ is an arbitrary point in $B_1$, we see that $H\equiv 0$ in $B_1$, a contradiction.

This completes the proof of theorem \ref{thm} \endproof

A sketch of proof of theorem \ref{thm2}. If $M$ is a  submanifolds in a space form $N(c)$, then it is biharmonic if and only if
\begin{eqnarray}
\Delta^\bot \vec{H}-\sum_{i=1}^mh(A_{\vec{H}}e_i, e_i)+cm\vec{H}&=&0, \label{equ4}
\\m\nabla|\vec{H}|^2+4\sum_{i=1}^mA_{\nabla^\bot_{e_i} \vec{H}}e_i&=&0. \label{equ5}
\end{eqnarray}

The same as the step 1 in the proof of theorem \ref{thm}, by using equation (\ref{equ5}), we can prove that if  $H(p)\neq 0$ at a point $p\in M$, then it is a constant around a neighborhood of the point $p$. Thus we have at $p$, $\Delta|\vec{H}|^2=0.$

On the other hand at $p$
\begin{eqnarray*}
\Delta|\vec{H}|^2&=&2|\nabla \vec{H}|^2+2\langle\vec{H}, \Delta\vec{H}\rangle
\\&=&2|\nabla \vec{H}|^2-2\sum_{i=1}^mR^N(e_i, \vec{H}, e_i, \vec{H})
\\&=&2|\nabla \vec{H}|^2-2c m|\vec{H}|^2
\\&\geq&2|\nabla \vec{H}|^2.
\end{eqnarray*}
Therefore  $\nabla\vec{H}(p)=0$. Now we choose an orthogonal basis $\{e_i, i=1,...,m\}$ of $T_pM$. Computing directly one gets at $p$
\begin{eqnarray}
0=\langle\nabla_{e_i}\vec{H}, e_j\rangle=H\langle\nabla_{e_i}\nu, e_j\rangle=-H\langle\nu, \nabla_{e_i}e_j\rangle=H\langle h(e_i, e_j),\nu\rangle,
\end{eqnarray}
for any $1\leq i,j\leq m$. Taking trace over this equality we get
$$H^2(p)=0.$$
Therefore $H(p)=0$. This is a contradiction to $H(p)\neq 0$.

This completes the proof of theorem \ref{thm2}. \endproof

\textbf{Acknowledgement.} Special thanks due to my supervisor professor Guofang Wang for his constant encouragement and support.
 {}
\vspace{1cm}\sc
Yong Luo

Mathematisches Institut,

Albert-Ludwigs-Universit\"at Freiburg,

Eckerstr. 1, 79104 Freiburg, Germany.

{\tt yong.luo@math.uni-freiburg.de}

\vspace{1cm}\sc

\end{document}